\documentstyle[amsfonts,amscd,amssymb]{amsart}
\newtheorem{thm}{Theorem}[section]
\newtheorem{lem}[thm]{Lemma}
\newtheorem{prop}[thm]{Proposition}
\newtheorem{cor}[thm]{Corollary}
\theoremstyle{definition}
\newtheorem{defn}{Definition}[section]
\newtheorem{ex}[defn]{Example}

\newtheorem{nt}[defn]{Notation}

\numberwithin{equation}{section}

\begin{document}
\newcommand{\Soc}{\operatorname{Soc}}
\newcommand{\modP}{\mod_{P}\Lambda}
\newcommand{\modl}{\mathrm{mod}\La}
\newcommand{\module}{\mod}
\newcommand{\Module}{\operatorname{Mod}}
\newcommand{\Cok}{\operatorname{Coker}}
\newcommand{\Hom}{\operatorname{Hom}}
\newcommand{\Ext}{\operatorname{Ext}}
\newcommand{\Tr}{\operatorname{Tr}}
\newcommand{\rad}{\operatorname {{\bold r}}}
\newcommand{\La}{\operatorname{\Lambda}}
\newcommand\End{\operatorname{End}}
\newcommand\ann{\operatorname{ann}}
\newcommand\Img{\operatorname{Im}}
\newcommand\D{\operatorname{D}}
\newcommand\Ker{\operatorname{Ker}}
\newcommand\Coker{\operatorname{Coker}}
\newcommand{\bt}{\operatorname {{\bold t}}}
\newcommand{\Gr}{\operatorname{Gr}}
\newcommand{\prd}{\operatorname{pd}}
\newcommand{\gd}{\operatorname{{gl.dim}}}
\newcommand{\radic}{\operatorname{rad}}
\newcommand{\ind}{\operatorname{ind}}
\newcommand{\cx}{\operatorname{c}}
\bibliographystyle{plain}

\title[Algebras with smallest resolutions]{Finite-dimensional algebras with smallest resolutions
of simple modules}
\author{Shashidhar Jagadeeshan}
\address{Centre for Learning\\ 462, 9th Cross \\ Jayanagar 1st
Block\\Bangalore 560 012 \\ India} \email{jshashidhar@@gmail.com}
\author{Mark Kleiner}
\address{Department of Mathematics, Syracuse University, Syracuse,
  New York 13244-1150}
  \thanks{The first author would like to thank the Department of
  Mathematics at Syracuse University for
  having given him the opportunity to spend the academic year
  2002-2003 there.}
\email{mkleiner@@sound.syr.edu}  \maketitle

\section*{Introduction}

Let $\La$ be an associative ring with identity and with the
Jacobson radical $\rad$,  let $\modl$ be the
category of finitely generated left $\La$-modules, and let
$\Lambda^{op}$ be the opposite ring of $\Lambda$. All modules are
left unital modules, and if $X$ is a module then $\prd X$ is the projective dimension of $X$.  If $\La$ is left artinian and $M\in\mathrm{mod}\La$, we denote by $P(M)$ a projective cover of $M$. Throughout the paper we fix an arbitrary
field $k$.   For terminology and notation, we refer the reader to
~\cite{ars,be91}.

If $\La$ is a finite-dimensional $k$-algebra, any nonzero
$M\in\text{mod}\La$ has a minimal projective resolution
\begin{equation}\label{resint}
\dots\to P_2 \overset{d_2}\to P_1\overset{d_1}\to P_0
\overset{d_0}\to M\to0,
\end{equation}
which is a very important but rather complicated homological
invariant of the  module $M$.  An interesting numerical invariant
of the above resolution is the complexity of $M$, $\cx_{\La}(M)$,
which  measures the growth of the size of
the $n$-th term, $P_n$, as $n\to\infty$.  The notion of complexity
is especially important when $\La=kG$ is the group algebra over
$k$ of a finite group $G$, where char$\,k$ divides the order of
$G$, because $\cx_{kG}(M)=\cx_{G}(M)$ is the dimension of the cohomological support
variety $V_{kG}(M)=V_G(M)$ associated to $M$ (see ~\cite[Ch.
5]{be91}).  We recall the definition of complexity.  Let $\mathbb
N$ be the set of nonnegative integers, and let ${\mathbf
b}=\{b_n\,|\, n\in\mathbb N\}$ be an arbitrary sequence with
$b_n\in\mathbb N$. The growth of $\mathbf b$ is the least
nonnegative integer $\gamma(\mathbf b)$ for which there is a
positive real number $A$ satisfying $b_n\le An^{\gamma({\mathbf
b})-1}$ for sufficiently large $n$.  If no such integer exists
then $\gamma(\mathbf b)=\infty$.  By definition,
$\cx_{\La}(M)=\gamma(\{\dim_k P_n\})$.

 For an
arbitrary finite-dimensional algebra $\La$,  the cohomological support variety $V_{\La}(M)$
of $M\in\text{mod}\La$ is defined in ~\cite{snsol04}.  The
definition works well for some classes of algebras; in particular,
many properties of the support variety over a group algebra extend
to self-injective algebras over an algebraically closed field that
satisfy appropriate finite generation assumptions ~\cite{ehsst}.
However, some of those properties do not extend to all
algebras.  For example, let $\La$ be a self-injective algebra whose Hochschild cohomology ring $\mathrm{HH}^*(\La)$ contains a commutative noetherian graded subalgebra $H$ where (i) $H^0=\mathrm{HH}^0(\La)=Z(\La)$ and (ii) $\Ext^*_{\La}(\La/{\mathbf r}, \La/{\mathbf r})$ is a finitely generated $H$-module.  Under these assumptions, if $M\in\text{mod}\La$ then $\cx_{\La}(M)<\infty$, and $V_{\La}(M)$ is trivial if and only if $M$ is projective. On the other hand, the authors of   ~\cite{snsol04} and ~\cite{ehsst} have constructed  finite-dimensional algebras $\La$ of infinite global dimension over an algebraically closed field for which  $\cx_{\La}(M)=\infty$ for some $M\in\text{mod}\La$, or $V_{\La}(M)$ is trivial for all $M\in\text{mod}\La$.  Some of these algebras belong to the class of algebras arising in this paper, and Example \ref{quiver8} exhibits two such algebras for which each module is eventually periodic (see the definition below).  However, for one of the algebras the support variety of each module is trivial, while for the other some modules have nontrivial support varieties. Thus, at least  for the algebras in this class, there is no clear connection between the cohomological support variety of a module and the behavior of the terms of its minimal projective resolution.  We also note that the notions of complexity and  support variety do not
distinguish between modules of finite projective dimension: all
such modules have  complexity zero and trivial support variety.  

To
study modules of arbitrary projective dimension over a class of finite-dimensional algebras containing all algebras over an
algebraically closed field, we examine a numerical invariant of
the resolution (\ref{resint}) that is finer than $\cx_{\La} (M)$.
This invariant is the sequence
\newline\centerline{$p(M)=\{p(M)_n=\ell(P_n)\,|\,n\in\mathbb
N\}$} \newline\noindent where $\ell(X)$ is the length of
$X\in\text{mod}\La$.

The obvious question is why $p(M)$ is preferable to $\{\dim_k
P_n\}$.  There are several answers. First of all, one can recover
$\cx_{\La}(M)$ from $p(M)$ because, according to
~\cite{av89a}, $\gamma(\{ \dim_k P_n\})=\gamma(p(M))$ so that
$\cx_{\La}(M)=\gamma(p(M))$, and the latter formula also works when
$\La$ is a left artinian ring. Second, $p(M)$ is a categorical
invariant of $M$, and thus is preserved under Morita equivalence,
while $\{\dim_k P_n\}$ is not.  Third, it may be easier to detect
finite projective dimension from $p(M)$ than from $\{\dim_k
P_n\}$, e.g., if $\ell(P_n)=1$ then $P_n$ is simple projective, so
$P_{n+t}=0$ for $t>0$, while it is impossible to detect a simple
projective module $P$ from $\dim_k P$ when $\dim_k P>1$.  Finally,
for  a minimal projective resolution
\begin{equation}\label{resint.5}
\dots\to Q_2 \to Q_1\to Q_0
\to N\to0
\end{equation} of $N\in\text{mod}\,\Gamma$ where $\Gamma$ is
a left artinian ring, we can compare the lengths $\ell(P_m)$ and
$\ell(Q_n)$, while the comparison of dimensions  of $P_m$ and
$Q_n$ over a field may be either impossible or unnatural.

Of course, certain properties of the sequence $p(M)$ have already
been addressed in the literature; see, for instance, ~\cite[Theorem 5.10.4 and Corollary 5.10.7]{be91} on the periodicity of $M\in\text{mod}\La$ when $\La=kG$ is the
group algebra and $\cx_{\La}(M)=1$. We recall that given the minimal  projective resolution
(\ref{resint}), the $n$-th syzygy of $M$ is $\Omega^n
M=\Coker d_{n+1}$, and $M$ is said to be {\it eventually periodic} if there exists an integer $t>0$  such that
$\Omega^{s+t} M\simeq\Omega^s M$ for some integer $s\ge 0$; if $s=0$ then $M$ is {\it periodic}. The smallest possible $t$ is the {\it period} of $M$.  

Our approach to the set of sequences $p(X)$ is based on the
introduction of a means to compare two such sequences, $p(X)$ and
$p(Y)$: we define  a preorder $\preccurlyeq$, i.e., a
reflexive and transitive binary relation, on the
set, and say that the minimal projective resolution of $X$ is
\lq\lq smaller" than that of $Y$ if $p(X)\preccurlyeq p(Y)$. There
are several ways to define the preorder $\preccurlyeq$, but all
preorders in this paper are such that if $p(X)_n\le p(Y)_n$ for
all $n\in\mathbb N$,  then $p(X)\preccurlyeq p(Y)$, i.e., we
respect the intuitive notion of what it means that $p(X)$ is
smaller than $p(Y)$.

The following easily verifiable statement shows that among the modules of fixed length, semisimple modules have the largest resolutions. 

\begin{prop}\label{intr1} Let $\La$ be a left
artinian ring, let $L,M,N\in\mathrm{mod}\La$, and let $n\in\mathbb
N$.
\begin{itemize}
\item[(a)] If $\,0\to L\to M\to N\to0$  is an exact sequence then
$p(M)_n\le p(L)_n+p(N)_n$. \item[(b)]  $p(M\oplus
N)_n=p(M)_n+p(N)_n$. \item[(c)] Let $S_1,\dots,S_t$ be a complete
set of pairwise nonisomorphic simple $\La$-modules
and let $m_i$ be the multiplicity of $S_i$ in a composition series of $M, 1\le i\le t$. Then
\newline$p(M)_n\le p(S_1^{m_1}\oplus\dots\oplus S_t^{m_t})_n
=m_1p(S_1)_n+\dots+m_t p(S_t)_n$.
\end{itemize}
\end{prop}

We are interested in those left artinian rings for which the largest projective resolutions  are as small as possible.  In view of Proposition \ref{intr1},  semisimple modules have the largest resolutions, so we wish to know when the minimal projective resolutions of their indecomposable summands, the simple modules, are as small as possible. More precisely, the problem is to  describe the left artinian rings $\La$ with the following property.\vskip.05in

\noindent$(*)$ \centerline{If $S$ is a simple $\La$-module and
$T\in\mathrm{mod}\,\Gamma$, where $\Gamma$ is any left
artinian ring,}  \centerline{then $\prd S\le\prd T$  implies $p(S)\preccurlyeq
p(T)$.}\vskip.05in

For a
preorder $\preccurlyeq$ described below, we solve the problem when
$\La$ is a  finite-dimensional $k$-algebra that  is  {\it
elementary}, which means that $\La/\rad$ is isomorphic to a direct
product of several copies of $k$.  This includes all  algebras
over an  algebraically closed field because each such algebra  is
Morita equivalent to an elementary algebra. Since a
finite-dimensional algebra is elementary if and only if it is
isomorphic to the path algebra of a finite quiver with relations
(see ~\cite[Section III.1]{ars}), we characterize the above
algebras in terms of quivers and relations (Theorem
\ref{quiver5}).  This is our main result.

We now describe the preorders $\preccurlyeq$ considered in this paper. For a left artinian ring
$\La$, we do not know which sequences ${\mathbf a}=\{a_n\,|\,
n\in\mathbb N\}$ of nonnegative integers satisfy  ${\mathbf a}=p(M)$ for some  nonzero $M\in\mathrm{mod}\La$.  So
we define a
preorder $\preccurlyeq$ on the collection $\mathcal{A}$ of all sequences
$\mathbf{a}$ that satisfy some obvious necessary conditions for
the existence of the indicated module $M$.  In fact, we define not
one, but infinitely many preorders that depend on the positive
integer parameter $r$.  When restricted to the sequences
associated with modules, the preorder works as follows. Given an
$r>0$ and $M,N\in\mathrm{mod}\La$, we have $p(M)\preccurlyeq p(N)$
if and only if $\prd M\le \prd N$ and, for all $n\in\mathbb N$, if
$p(M)_n> p(N)_n$ then $\prd M<n+r$ {(Definition \ref{res3})}.
Note that if $\prd M<r$ and $\prd N\ge r$, then $p(M)\preccurlyeq
p(N)$ but $p(N)\not\preccurlyeq p(M)$, i.e., the preorder
distinguishes between some modules of finite projective dimension, and if  $\prd M= \prd N=\infty$ then
$p(M)\preccurlyeq p(N)$ if and only  if, for all $n\in\mathbb N$,
$p(M)_n\le p(N)_n$. In particular, the preorder (for any $r>0$) is easy to understand if $\La$ is a selfinjective $k$-algebra:  $p(M)\preccurlyeq p(N)$ whenever $M$ is projective, and if $M$ is not projective then $p(M)\preccurlyeq p(N)$ if and only  if, for all $n\in\mathbb N$, $p(M)_n\le p(N)_n$ . 

For the preorders just described, a left artinian ring $\La$ satisfies $(*)$ if and only if it satisfies the following seemingly weaker condition. \vskip.05in

\noindent$(**)$ \centerline{If $S$ is a simple $\La$-module and
$T$ is a simple $\Gamma$-module, where $\Gamma$ is any left
artinian ring,}  \centerline{ then $\prd S\le\prd T$  implies $p(S)\preccurlyeq
p(T)$.}\vskip.05in

\noindent This equivalence of  $(*)$ and  $(**)$ follows from a careful examination of the proofs in Section 1, we have no a priori explanation for it.

It is easy to describe the elementary $k$-algebras satisfying $(*)$ when $r=1$, our main result concerns the value $r=2$, and the case $r>2$ is open.

For  $r=2$ and $\La$  a left
artinian ring, condition $(*)$  is equivalent (Corollary \ref{res6}(iii)) to the following.\vskip.05in

\noindent (F) \centerline{The radical of each indecomposable projective
module is either projective or simple.} 
\vskip.05in
\noindent Condition (F) is similar to other conditions that have appeared in the literature.  
One such condition is the following.\vskip.05in

\noindent (A) \centerline{Each submodule of an indecomposable projective
module is either projective or simple.} \vskip.05in

\noindent The structure of left artinian rings satisfying (A) is described in  ~\cite{jk2} as part of the description of projectively stable rings, i.e., rings for which a morphism of finitely generated modules without nonzero projective direct summands must be zero if it factors through a projective module. Consider the following condition.\vskip.05in

\noindent (B) \centerline{The injective envelope of each simple nonprojective torsionless module is projective. } \vskip.05in

\noindent Then a left artinian ring is projectively stable if and only if it satisfies (A) and (B).
Another condition similar to (F) is the following.\vskip.05in

\noindent (D) \centerline{Each indecomposable submodule of an indecomposable projective module is} \centerline{either projective or simple. }\vskip.05in

\noindent This condition has to do with rings that are stably equivalent to left hereditary rings. Consider the following condition.\vskip.05in

\noindent (E) \centerline{Each simple nonprojective torsionless module is cotorsionless. }

\vskip.05in 

\noindent Then a left artinian ring is  stably equivalent to a left hereditary left artinian ring if and only if it satisfies (D) and (E); this is presented in  ~\cite[p.40]{ar} for artin algebras, but is known to be true in the indicated setting.

Clearly, (A)$\Rightarrow$(F) but (F)$\not\Rightarrow$(A), (D)$\not\Rightarrow$(F), and (E)$\not\Rightarrow$(F).  We prove that (F)$\Rightarrow$(D) (Proposition \ref{proj1}).  However, a finite-dimensional algebra satisfying (F) need not be stably equivalent to a hereditary algebra, as follows from our main result and ~\cite{br}.  Hence (F)$\not\Rightarrow$(E).  

The main result of this paper says that an elementary $k$-algebra satisfies (F) if and only if it is isomorphic to a quadratic monomial algebra (see the definition in Section 3) determined by a (unique) triple $(G,\mathcal{D},\rho)$ where  $G$ is a finite quiver (directed graph) with the set of vertices $v(G)$ and without oriented cycles, $\mathcal{D}$ is a subset of $v(G)$ consisting of sinks (vertices at which no arrow starts), and $\rho:{\mathcal D}\to v(G)$ is a function such that the sinks in  $\Img \rho$ belong to $\mathcal{D}$.  A simple geometric construction produces the quadratic monomial algebra from $(G,\mathcal{D},\rho)$.

Section 1 of the paper defines the set
of sequences $\mathcal{A}$  and preorders
$\preccurlyeq$ as indicated above, and we show for all of these preorders that if $\La$ satisfies $(*)$, then every semisimple module in $\mathrm{mod}\La$ is eventually periodic, and every $M\in\mathrm{mod}\La$ satisfies $\cx_{\La}(M)=1$.  The section also contains a 
computation (Corollary \ref{res8}) of the global dimension of a left artinian ring
satisfying (F) in terms of sequences of morphisms between
indecomposable projective modules of length two.  In Section 2, for a left artinian ring $\La$ satisfying (F), we show that every $M\in\mathrm{mod}\La$ is eventually periodic, and develop
properties of projective $\La$-modules needed for the proof of the main result in Section 3.   Condition (F) is not selfdual in that an algebra $\La$ may
satisfy (F), while $\La^{op}$ does not (Example \ref{quiver8}(c)), so we characterize (Proposition \ref{quiver7})
all elementary $k$-algebras $\La$ such that both $\La$ and  $\La^{op}$
satisfy (F). The point of Section 4 is that an elementary $k$-algebra $\La$ satisfying (F) not only has well-behaved projective resolutions of its modules,  but possesses other good homological properties.  Since such an algebra $\La$ is quadratic monomial, it is a Koszul algebra, and its Ext-algebra
$E(\La)={\overset{\infty}{\underset{i=0}{\oplus}}}\Ext^{i}_{\La}(\La/\rad,
\La/\rad)$ is finitely generated ~\cite[Proposition 2.2]{gz}.  We show that $E(\La)$ is left noetherian, which happens rarely, so considerations of  ~\cite[Section 4]{mz} apply.  In particular, the Poincar\'e series of a finitely  generated graded $\La$-module is rational.  We note that although the authors of  ~\cite{mz} assume $E(\La)$ noetherian, their proof only uses that $E(\La)$ is left noetherian.

We are grateful to the referee for the very helpful suggestions.  In particular, the referee asked the following very natural question.  For which elementary algebras $\La$ with property (F) does the Koszul dual $E(\La)$ share this property?  The question is suggested by our characterization of those algebras with property (F) for which the opposite algebra shares this property.  Since $E(\La)$ is in general infinite-dimensional (resp. not left artinian), our results do not apply immediately because we always assume finite-dimensionality (resp. left artinianness).  However, the referee was wondering to what extent the methods and constructions of this paper actually depend on these assumptions.  Some few small examples the referee looked at do show that Koszul duals can very well share property (F).  We intend to investigate this question in the future.

\section{Sequences and projective resolutions}\label{sequences}

In this section we fix a left artinian ring $\La$.

Let $\mathcal{A}$ be the set of infinite sequences ${\mathbf a} = \{ a_n
\,|\,n\in{\mathbb N}\}$ of nonnegative integers with the property that $a_0 > 0$
and, for all $n$, if $a_n \le 1$ then $a_i =0$ for $i> n$. We set
$ \dim {\mathbf a} = \sup \, \{ n \, |\, a_n > 0 \}$, ${\mathcal{A}}_n=\{ {\mathbf a} \in {\mathcal A}|\dim {\mathbf a}
\ge n \}$ for $0 \le n <
\infty$, and ${\mathcal{A}}_{\infty}$ = ${\underset{n =0}{
\overset{\infty}{\bigcap}}}$ ${\mathcal{A}}_{n}$. We have
$0 \le \dim {\mathbf a} \le \infty $; $\dim {\mathbf a} =n<\infty$ if and only if ${\mathbf a}
\in$ ${{\mathcal{A}}_n} - {{\mathcal{A}}_{n+1}}$; and $\dim {\mathbf a} =
\infty$ if and only if ${\mathbf a} \in$ ${\mathcal{A}}_{\infty}$.
Our motivation for studying such sequences comes from the obvious fact that if 
$ 0 \neq M\in \mathrm{mod}\La$, then $ p(M) = \{\ell (P_{n}) \} \in
\mathcal{A}$ and $\prd
M =\dim p(M)$.

In order to introduce appropriate preorders on
${\mathcal{A}}$,  we need the following easy statement.  For an arbitrary preorder $(\mathcal{P},\preccurlyeq)$ we say that an element $x\in\mathcal{P}$ is {\it a least} element if $x\preccurlyeq p$ for all $p\in\mathcal{P}$; if $x,y\in\mathcal{P}$ are least elements, then $x\preccurlyeq y$ and $y\preccurlyeq x$.

\begin{lem}\label{res1}
\begin{itemize}
     \item [(a)]  Let $0\le m<\infty$ be an integer. Denote by $\Gamma_1$ the quotient of the path algebra over $k$ of the quiver with  vertices $v_1,\dots,v_{m+1}$ and arrows $\alpha_i:v_i\to v_{i+1},\ i=1,\dots,m$, modulo the ideal generated by the elements $\alpha_{i+1}\alpha_i,\ i=1,\dots,m-1$. If $T_1$ is the simple $\Gamma_1$-module associated with the vertex $v_1$, then $p(T_1)={\mathbf t_1}=\{a_n\}$ where $a_n=2$ for $n<m$, $a_m=1$, and $a_n=0$ for $n>m$.  In particular, $\prd T_1=m$.
     \item [(b)] Let $\Gamma_2=k[X]/(X^2)$ where $k[X]$ is the polynomial algebra. If $T_2$ is the simple $\Gamma_2$-module, then $p(T_2)={\mathbf t_2}=\{b_n\}$ where $b_n=2$ for all $n$.  In particular, $\prd T_2=\infty$. 
	\item [(c)] Let $\preccurlyeq$  be a preorder on ${\mathcal{A}}$ such that ${\mathbf x}\preccurlyeq{\mathbf y}$ whenever $x_n\le y_n$ for all $n$.  For a fixed $0\le m<\infty$, ${\mathbf t_1}$ is a least element of ${\mathcal{A}}_m$.  The element ${\mathbf t_2}$ is a least element of ${\mathcal{A}}_\infty$. 
\end{itemize}
\end{lem}
\begin{pf} (a) and (b)  We leave the easy verification to the reader.

(c) If ${\mathbf c}\in{\mathcal{A}}_m$, then $c_m>0$ whence $c_n>1$ for $n<m$.  Then $a_n\le c_n$ for all $n$ and we conclude that  ${\mathbf t_1}\preccurlyeq{\mathbf c}$.  If ${\mathbf c}\in{\mathcal{A}}_\infty$, then  $c_n>1$ for all $n$, whence $b_n\le c_n$ for all $n$, so  ${\mathbf t_2}\preccurlyeq{\mathbf c}$.
\end{pf}	

We get the most simple-minded preorder on $\mathcal{A}$ by setting ${\mathbf a} \le {\mathbf b}$ if and only if $a_n
\le b_n$ for all $n$. Then  $\le$ is a partial order and we have
the following statement.
\begin{prop}\label{res2} For the left artinian ring  $\La$ and partial order $(\mathcal{A}, \le)$, the following are equivalent. 
\begin{itemize}
     \item [(a)] For each simple $\La$-module $S$, $p(S)$ is the least element of ${\mathcal{A}}_{\prd S}$.
	\item [(b)] $\La$ satisfies $\mathrm{(*)}$.   
	\item [(c)]  For each
indecomposable projective $\La$-module $P$, $\ell(P) \le 2$.
\end{itemize} 
\end{prop}
\begin{pf} (a)$\Rightarrow$(b) Obvious.

(b)$\Rightarrow$(c) Let $\La$ satisfy $\mathrm{(*)}$ and let $S$ be a simple $\La$-module with $m=\prd S$.  If $m<\infty$ then $p(S)\le p(T_1)={\mathbf t_1}$ where $T_1$ is defined in Lemma \ref{res1}(a), whence $\ell(P(S)) \le \ell(P(T_1))=2$.  If $m=\infty$ then $p(S)\le p(T_2)={\mathbf t_2}$ where $T_2$ is defined in Lemma \ref{res1}(b), whence $\ell(P(S)) \le \ell(P(T_2))=2$.

(c)$\Rightarrow$(a)  Let $S$ be a simple $\La$-module with $m=\prd S$.  Since $\ell(P) \le 2$ for each indecomposable projective $\La$-module $P$, an easy inductive argument on the minimal projective resolution of $S$ shows that $p(S)=p(T_1)={\mathbf t_1}$ if $m<\infty$, and $p(S)=p(T_2)={\mathbf t_2}$ if $m=\infty$.  By Lemma \ref{res1}(c), $p(S)$ is a least element of ${\mathcal A}_m$.  Since $({\mathcal A},\le)$ is a partial order, $p(S)$ is the least element.
\end{pf}

Proposition \ref{res2} shows that for the partial order $(\mathcal{A},\le)$, the left artinian rings satisfying $(*)$ form a rather small class of
radical-square-zero rings. To produce more interesting classes of rings, we need coarser preorders on $\mathcal{A}$.

\begin{defn}\label{res3}
Let $r>0$ be a fixed integer. For ${\mathbf a, b} \in \mathcal{A}$
we set ${\mathbf a} \preccurlyeq  {\mathbf b}$ if and only if $\dim {\mathbf a}\le\dim {\mathbf b}$ and, for all $n$, $a_n>  b_n$ implies $\dim {\mathbf a}<{n+r}$.
\end{defn}

Applying Definition \ref{res3} to the minimal projective resolutions (\ref{resint}) and  (\ref{resint.5}) of modules $M$ and $N$,  we see that ${p(M)}\preccurlyeq\ {p(N)}$ if and only if  (i) $\prd M\le\prd N$ and (ii) the violation of the intuitive condition $\ell(P_n)\le \ell(Q_n)$ for some $n$ is allowed only if $\prd M<\infty$, and for at most $r$ terms at the tail of the resolution of $M$.  Note also that if we formally apply Definition \ref{res3} to $r =0$, we will obtain
the partial order $(\mathcal{A},\le)$. 

\begin{prop}\label{res4} Let $r>0$ be a fixed integer, let
$(\mathcal{A},\preccurlyeq)$  be given by Definition \ref{res3},
and let ${\mathbf a, b} \in \mathcal{A}$. 
\begin{itemize}
     \item [(a)] The binary relation $(\mathcal{A},\preccurlyeq)$ is a preorder.
     \item [(b)] Both  ${\mathbf a} \preccurlyeq {\mathbf b}$ and ${\mathbf b} \preccurlyeq {\mathbf a}$ hold if and only if $ \dim {\mathbf a} = \dim {\mathbf b}$ and $a_i = b_i$ for $ 0 \le i \le (\dim
     {\mathbf a}) - r$.
     \item [(c)] For $ 0 \le m < \infty, {\mathbf b} \in {\mathcal{A}}_{m}$
     is a least element of ${\mathcal{A}}_m$ if and only if $\dim {\mathbf b} =m$ and
     $b_{i}=2$ for $ 0 \le i \le m-r$.  An element ${\mathbf b} \in{{\mathcal{A}}_{\infty}}$ is a least element
     of ${{\mathcal{A}}_{\infty}}$ if and only if
$b_{i} =2$ for $ 0 \le i < \infty$.
\end{itemize}
\end{prop}
\begin{pf} (a) The reflexivity of $\preccurlyeq$ is obvious. To
prove the transitivity, let ${\mathbf a} \preccurlyeq {\mathbf b}$  and ${\mathbf b}
\preccurlyeq {\mathbf c} $.  Then $\dim{\mathbf a}\le\dim{\mathbf b}$ and $\dim{\mathbf b}\le\dim{\mathbf c}$ whence $\dim{\mathbf a}\le\dim{\mathbf c}$. Let $ a_n > c_n$ for some $
n \ge 0$. If $ b_n \le c_n$ then $a_n > b_n$, whence $\dim{\mathbf a}<n+r$ because ${\mathbf a} \preccurlyeq {\mathbf b}$.  If  $ b_n > c_n$, then $ \dim{\mathbf b}<n+r$ because ${\mathbf b}
\preccurlyeq {\mathbf c} $ so $\dim{\mathbf a}<n+r$.  Since $\dim{\mathbf a}<n+r$, we get ${\mathbf a} \preccurlyeq {\mathbf c}$.

(b) This is clear.

(c)  Let $m<\infty$. If ${\mathbf  b}\in {\mathcal{A}}_{m}$ is a least element of ${\mathcal{A}}_{m}$, we must have ${\mathbf  b}\preccurlyeq{\mathbf  t_1}$ where ${\mathbf  t_1}$ is defined in Lemma \ref{res1}(a). By Lemma \ref{res1}(c), we also have ${\mathbf  t_1}\preccurlyeq{\mathbf b}$. Now the statement is an immediate consequence of  (b) and the  definition of ${\mathbf  t_1}$.  The same argument with the replacement of ${\mathbf  t_1}$ by ${\mathbf  t_2}$ works in the case ${\mathbf b} \in{{\mathcal{A}}_{\infty}}$.
\end{pf}

The following statement extends Proposition \ref{res2}.  

\begin{prop}\label{res5} Let $r>0$ be a fixed integer. For the left artinian ring  $\La$ and preorder $(\mathcal{A},\preccurlyeq)$, the following are equivalent.
\begin{itemize}
     \item [(a)] For each simple $\La$-module $S$, $p(S)$ is a least element of
${\mathcal{A}}_{\prd S}$.
     \item [(b)] The ring $\La$ satisfies $(*)$.

\item [(c)]  For each simple $\La$-module $S$, either $\prd S < r$ or $\ell(P(S)) =2$.
\end{itemize}
\end{prop}
\begin{pf} (a) $\Rightarrow$ (b) Obvious.

(b) $\Rightarrow$ (c) Let $S$ be a simple module with $m=\prd S \ge r $.   If $m<\infty$ then $p(S)\preccurlyeq p(T_1)$ where $T_1$ is defined in Lemma \ref{res1}(a). Since $m\ge r$, we cannot have $p(S)_0>p(T_1)_0$. Hence $\ell(P(S))=p(S)_0\le p(T_1)_0=2$.  If $m=\infty$ then $p(S)\preccurlyeq p(T_2)$ where $T_2$ is defined in Lemma \ref{res1}(b), whence $\ell(P(S)) \le\ell(P(T_2))= 2$.

(c) $\Rightarrow$ (a) Let  $S$ be a simple $\La$-module, and let $\prd S=m\le \dim{\mathbf a}$ where ${\mathbf a}\in {\mathcal{A}}$. If $m<r$, then $m-r<0$ so no integer $i$ satisfies $0\le i\le m-r$. If $m\ge r$ then $p(S)_0=\ell(P(S))=2$, and we obtain by induction on $i$ that if $m<\infty$ then $p(S)_i=2$ for $0\le i\le m-r$, and  if $m=\infty$ then $p(S)_i=2$ for $0\le i<\infty$. By Proposition \ref{res4}(c), $p(S)$ is a least element of ${\mathcal{A}}_{\prd S}$.
\end{pf}

\begin{cor}\label{res5.5} For a fixed integer $r>0$, denote by $\preccurlyeq$ the preorder defined in Definition \ref{res3} that depends on $r$, and suppose that $\La$ satisfies $(*)$ relative to $\preccurlyeq$.
\begin{itemize}
\item [(a)] Every semisimple module in $\mathrm{mod}\La$ is eventually periodic.
 \item[(b)] For all $M\in{\mathrm{mod}}\La$, $\cx_{\La}(M)=1$.
\end{itemize}
\end{cor}
\begin{pf}  (a) It suffices to show that every simple $\La$-module is eventually periodic, which follows from Proposition \ref{res5}(c): if $\prd S=\infty$, every projective in the minimal projective resolution of $S$ is of length $2$, so all syzygies are simple.

(b) This is a consequence of (a) and Proposition \ref{intr1}(c).
\end{pf}
We write $\gd \La$ for the global dimension of $\La$.

\begin{cor}\label{res6} Let $r>0$ be a fixed integer.
\begin{itemize}
\item [(i)] If  $\gd \La < r$ then $\La$
satisfies $(*)$.
 \item[(ii)] If
$r=1$, then $\La $ satisfies $(*)$ if and only if $\ell(P) \le 2$ for each indecomposable
projective $\La$-module $P$.
  \item [(iii)] If $r=2$, then $\La$ satisfies
  $(*)$ if and only if
it satisfies (F).
\end{itemize}
\end{cor}
 \begin{pf} The proof is a straight forward verification of condition (c) of Proposition \ref{res5}(c).  Note that for $r=1$, the class of left artinian
rings satisfying $(*)$ coincides with that described in Proposition \ref{res2}. 
 \end{pf}
 
For a fixed $r \ge 2$, we now  
estimate $\gd \La$ in the case when $\gd \La \ge r$ and $\La$ satisfies $(*)$.  We need to recall some definitions from ~\cite{jk1}. Let $\mathcal{Q}$ be the
full subcategory of $\text{mod} \La$ determined by the
indecomposable nonhereditary projective modules of length 2; recall that a projective module is hereditary if all its submodules are projective.  A  path  of length $n \ge 0$ in $\mathcal{Q}$ is a
sequence
\begin{equation}\label{res6.5}
 Q_n\overset{f_n}{\longrightarrow}Q_{n-1} \longrightarrow \dots
\longrightarrow Q_{1}\overset{f_1}{\longrightarrow}Q_0
 \end{equation} 
 where, for all $i$,
$Q_i\in\mathcal{Q}$ and $f_i$ is a nonzero nonisomorphism. If $\gd \La \ge r$ and $\La$ satisfies $(*)$, $\mathcal{Q}$ is not empty according to  Proposition \ref{res5}.
Set $ l = \sup \{\mathrm{lengths\ of\ paths\ in\ }\mathcal{Q} \} $.

\begin{prop}\label{res7} Let $r \ge 2$ be a fixed integer. If the  left artinian ring $\La$  satisfies $(*)$ and
has the property that $\gd \La \ge r$, then $ l +2 \le \gd \La \le l+r$.
\end{prop}
\begin{pf} If $\gd \La = m$, then $\prd S = m$ for some simple $\La$-module $S$. If $m=\infty$, a minimal projective resolution of $S$ is an infinite path in  $\mathcal{Q}$ in view  of Propositions \ref{res5} and \ref{res4}(c). Hence $l = \infty$ and the statement holds. 
If $r \le m  < \infty$, the same argument shows that the first $m-r+1$ terms of a minimal projective resolution of $S$ form a path in  $\mathcal{Q}$, whence $m-r\le l$ and $m\le l+r$. It is clear that any path in $\mathcal{Q}$ of the form (\ref{res6.5}) is an exact sequence of modules, whence $T= \Coker f_1$ and $\Ker f_n$ are simple nonprojective $\La$-modules, so $n+2\le\prd T\le m$.  Since $l=\sup\{n\}$, $l+2\le m$.
\end{pf}
\begin{cor}\label{res8} If $\La$ is nonhereditary and
satisfies (F), then $\gd \La = l + 2$.
\end{cor}
\begin{pf} We set $r=2$ and use Corollary \ref{res6} (iii) and Proposition \ref{res7}. Another way to prove the statement is to use ~\cite[Propositions
4.1 and 4.2]{jk1} and replace condition (A) by (F).
\end{pf} 

For the rest of the paper we focus on the left artinian rings satisfying
(F).

\section{Projective modules over $\La$ satisfying (F)}
Throughout this section $\La$ is a left artinian ring. We begin by rephrasing ~\cite[Ch. IV, Lemma 2.2]{ar}.
\begin{lem}\label{proj.5} Let $M$ be an indecomposable $\La$-submodule of $\oplus^n_{i=1}Q_i$ where each $Q_i$ is an indecomposable projective $\La$-module whose all submodules are either projective or simple.  Then $M$ is either projective or simple.
\end{lem}

\begin{prop}\label{proj1}
If $\La$ satisfies (F) then it satisfies (D).
\end{prop}
\begin{pf}  We show by induction on $\ell(P)$ that each indecomposable projective $\La$-module $P$ satisfies the requirement of (D).  This is clear when $\ell(P)=1$ or $\ell(P)=2$.  Let $\ell(P) > 2$ and let $M$ be an indecomposable nonprojective submodule of $P$.  Since $\ell(P) > 2$, $\rad P$ is not simple, hence is projective by (F), so $\rad P =
\oplus^n_{i=1}Q_i$ where each $Q_i$ is indecomposable projective.  Since $M$ is not projective, $M\subset\rad P$.  Since $\ell(Q_i) < \ell(P)$, the induction
hypothesis says that each indecomposable submodule of $Q_i$ is either projective or simple. By Lemma \ref{proj.5}, $M$ is simple.
\end{pf}

\begin{cor}\label{proj1.5} If $\La$ satisfies (F) then every $M\in\mathrm{mod}\La$ is eventually periodic.
\end{cor}
\begin{pf} By Proposition \ref{proj1}, $\Omega^1 M=P\oplus T$ where $P$ is projective and $T$ is semisimple.  Now the statement follows from Corollary \ref{res5.5}(a).
\end{pf}

Denote by ${\bold{a}}_{\La}$ the two-sided ideal of $\La$ equal to
the sum of nonprojective submodules  of $\Soc \La$. Recall that a module is torsionless if it is a submodule of a finitely generated projective module.

\begin{prop}\label{proj3} Suppose $\La$ satisfies (F).
\begin{itemize}
\item[(a)] The ring $\La / {\bold{a}}_{\La}$ is left hereditary. 
\end{itemize}
Let $S$
be a simple nonprojective  $\La$-module with the projective cover $P(S)$.
\begin{itemize}
\item[(b)] If $S$ is torsionless, then $S \simeq \rad P$ for some indecomposable projective $\La$-module $P$.
\item[(c)] The module $S$ is a projective $\La /{\bold{a}}_{\La}$-module if and only if $\ell(P(S))=2$ and $\Soc P(S) $ is  a
nonprojective $\La$-module.
\end{itemize}
\end{prop}

\begin{pf} (a) According to ~\cite[Proposition 5.3]{pl}, if $\La$ satisfies (D) then $\La / {\bold{a}}_{\La}$ is hereditary. It remains to use Proposition \ref {proj1}.

(b) Let $P$ be an indecomposable projective module of the smallest length
with $S \simeq T$ and $ T \subset P$.  Since $S$ is not projective, $ T \subset \rad P$.  If ${\rad}P$ is projective, $S \simeq U$ and $ U \subset Q$ for some indecomposable summand $Q$ of $\rad P$, which is impossible because  $\ell(Q) < \ell(P)$.  Thus ${\rad}P$ is not projective, so must be simple by (F), whence $S \simeq \rad P$. 

(c) The sufficiency is straight forward. For the necessity, we
note that \newline$P(S)/{\bold{a}}_{\La} P(S) \simeq S$ and
${\bold{a}}_{\La} P(S) \neq 0$.  Therefore $\rad P(S) = {\bold{a}}_{\La}
P(S)$ and, by ~\cite [Lemma 5.1] {pl}, ${\bold{a}}_{\La} P(S)$ is not
projective. Since $\La$ satisfies (F), $ {\bold{a}}_{\La} P(S) $ is
simple, so $\ell(P(S)) =2$.
\end{pf}

\section{Quivers and relations for algebras satisfying (F)}

In this section, $\La$ is an elementary finite-dimensional
$k$-algebra with the Jacobson radical $\rad$.
  
Denote by $H=(v(H), a(H))$ a finite
quiver with the set of   vertices $v(H)$ and the
set of arrows $a(H)$.  For each $\alpha\in a(H)$,
$s(\alpha)$ ($t(\alpha)$) is the starting (terminal) vertex of
$\alpha$. A vertex $x\in v(H)$ is a source (sink) if   no arrow
$\alpha\in a(H)$ satisfies $t(\alpha)=x$ ($s(\alpha)=x$).  If $p$ is a path in $H$,  $s(p)$ ($ t(p)$) is the starting (terminal) vertex of $p$, and we write $p:s(p)\to t(p)$.  For each $x\in v(H)$,
$e_x$ stands for the trivial path at $x$; a nontrivial path
consists of at least one arrow. 
The length of a path is the number of arrows in the path,
and a path $p$ is an oriented cycle if $p$ is  nontrivial and $s(p) =t(p)$. 

We recall (see ~\cite[Section III.1]{ars}) that a $k$-algebra
$\La$ is finite-dimensional and elementary if and only if it is isomorphic to the algebra
$k[H,\rho(H)]$ for some quiver $H$ with a
set of relations   $\rho(H)$.  Here $k[H,\rho(H)]$ is the quotient of the path
algebra $k H$ of $H$ modulo the two-sided ideal  $\langle\rho(H)\rangle$ generated by $\rho(H)$, where $\langle\rho(H)\rangle$ and the two-sided ideal $\langle a(H)\rangle$ generated by $a(H)$ satisfy $\langle
a(H)\rangle^m\subset\langle\rho(H)\rangle\subset \langle
a(H)\rangle^2$ for some positive integer $m$.  The image of  $a\in k H$ under the natural projection $\pi_{\rho}: k H
\to k[H,\rho(H)]$ is denoted by $\overline {a}$. 

We characterize
algebras satisfying (F) in terms of the quiver $H$ and relations
$\rho(H)$.  Of particular
interest to us are {\it quadratic monomial} algebras, for which
each element of $\rho(H)$ is a path of length 2. 

In the following definitions $\rho$ is a function whose domain is a set $\mathcal{D}$ of sinks of a quiver  $G$, and $\rho(H)$ is the set of relations on another quiver $H$ that is constructed from the triple $(G, \mathcal{D}, \rho)$.

\begin{defn}\label{def1}
Let  $(G, \mathcal{D}, \rho)$ be a triple consisting of a finite quiver $G$ without oriented cycles,  a subset $\mathcal{D}$ of $v(G)$
consisting of  sinks, and a  function $\rho:
{\mathcal{D}}\to v(G)$ for which all sinks in $\Img {\rho}$ belong to $ \mathcal{D}$.
A morphism $\varphi : (G, {\mathcal D},\rho) \to (G', \mathcal{D}',\rho')$ is a
morphism $\varphi : G \to G'$ of quivers satisfying $\varphi
(\mathcal{D}) \subset \mathcal{D'}$ and $\rho'\varphi|\mathcal{D}
= \varphi\rho$.
\end{defn}

\begin{defn}\label{def2}
The algebra $k(G,\mathcal{D}, \rho)$ of the triple $(G, \mathcal{D}, \rho)$ is the
algebra $k[H, \rho(H)]$, where the quiver $H$ and the set of relations $\rho(H)$ are obtained as
follows. For each $d \in \mathcal{D}$ we add a new single
arrow $\alpha_{d}: d\to
\rho(d)$, and then set $v(H) =v(G)$, $a(H) = a(G) \cup \{ \alpha_{d}\, |\, d
\in \mathcal{D} \} $,   and $\rho (H)=\{ \beta \alpha_{d}\,|\, \beta \in a(H),d \in\mathcal{D}\}$. Note that if $\mathcal{D}=\emptyset$ then $k(G,{\mathcal{D}}, \rho)=kG$.
\end{defn}

\begin{prop}\label{noloop} Let $ \La=k(G,{\mathcal{D}},\rho) = k[H, \rho(H)]$, $\rad={\mathrm{rad}}\La$, and $z\in v(H)$.
\begin{itemize}
   \item [(a)] Every oriented cycle in $H$ passes through a vertex
   in $\mathcal{D}$.
   \item [(b)]$\La$ is finite-dimensional.
   \item [(c)] If $z \in \mathcal{D}$ then $\rad {\bar{e}}_{z}$ is
   simple nonprojective.
   \item [(d)] If $ z \notin \mathcal{D}$ then $\rad {\bar{e}}_{z}$
   is projective.
\end{itemize}
\end{prop}
\begin{pf}
(a) Since $G$ has no oriented
cycles, an oriented cycle in $H$ must contain at least one of the added arrows $\alpha_{d}, \ d \in
\mathcal{D}$.

(b) Follows from (a) and the definition of $\rho(H)$.

Let  $W_z$ be the
set of those nontrivial paths in $H$ that start at $z$ and
have no subpath belonging to $\rho(H)$. Clearly, $ {\overline
{W}}_z = \{ \bar{p} \ | p \in W_z \}$ is a $k$-basis for $\rad {\bar{e}}_z={\mathrm {rad}}\La {\bar{e}}_z$.
For each $p \in W_z$ we have  $p = q \gamma $, where $\gamma \in
a(H)$  and $q$ is a path in $H$. 

(c) If $z=d \in \mathcal{D} $ and $p\in  W_z$ then, by construction,  $\gamma =
\alpha_{d}$ and $q =e_{\rho(d)}$. Hence $\rad {\bar{e}}_z=k{\bar{\alpha}}_d$ is 1-dimensional and, therefore, a simple $\La$-module. It is not
projective because $\rho(d)$ is not a sink in $H$.

(d) Suppose $ z \notin \mathcal{D} $.  If $z$ is a sink in $H$,  $\La {\bar{e}}_z$ is simple projective so $\rad {\bar{e}}_{z}=0$.  If $z$ is not a sink, let $\gamma_1, \dots,
\gamma_m,\ m>0$, be the arrows in $a(H)$ (and in $a(G)$) starting at $z$.  Then
$W_z = \overset{m}{\underset{j=1} {\bigcup}}  \{
q \gamma_j \  | q \in W_{t(\gamma_j)} \}  \cup \{\gamma_1,\dots,\gamma_m\}$, so $\rad {\bar{e}}_{z}
=\overset{m}{\underset{j=1} {\amalg}} \La {\bar{e}}_{t(\gamma_j)}
\ \bar{\gamma_j} \simeq \overset{m}{\underset{j=1} {\amalg}} \La
{\bar{e}}_{t(\gamma_j)} $  is projective.
\end{pf}

\begin{cor}\label{quiver1}
If $\La = k(G, \mathcal{D},\rho)$ then $\La$ satisfies (F).
\end{cor}
\begin{pf}
This follows from parts (c) and (d) of Proposition \ref{noloop} .
\end{pf}

Our goal now is to prove the converse of Corollary \ref{quiver1}.

\begin{nt} \label{not1} For the rest of the section, we fix a finite quiver $H$ with a set of
relations $\sigma(H)$, and assume that the algebra $\Sigma = k[H, \sigma(H)]$ is finite-dimensional and 
satisfies (F); the image of $a\in kH$ under the natural projection  $\pi_{\sigma}: kH \to \Sigma$ is denoted by $\hat{a}$. Let $\bold{s}={\mathrm {rad}}\,\Sigma$ and, for each $z \in
v(H)$, denote by $S_z$ (respectively, $P_z$)  the associated simple
(indecomposable projective) $\Sigma$-module.   We set
\begin{equation}\label{quiver1.1}
{\mathcal{D}}
= \{ d \in v(H) \ |\  \ell(P_d) =2,\  \Soc P_d {\mathrm{\ is\ not\
projective}}\}
\end{equation}
and define a function $\rho: {\mathcal{D}}\to v(H)$ by $\rho(d)
= y$ if $S_y \simeq \Soc P_d$.
\end{nt}

\begin{lem}\label{quiver2} In the setting of Notation \ref{not1}, we have:
\begin{itemize}
\item[(a)] For each $d \in \mathcal{D}$, there exists a unique
arrow $\alpha_{d}  \in a(H)$ such that $s(\alpha_{d}) =d$. We must
have $t(\alpha_{d}) = \rho(d)$. 
\item[(b)] If $d  \in
\mathcal{D}$ and $\beta \in a(H)$ satisfies $s(\beta)=\rho(d)$, then $\beta \alpha_d \in
\langle\sigma(H)\rangle$.

\end{itemize}
\end{lem}
\begin{pf} The statement is an easy consequence of $\ell(P_d) =2$.
\end{pf}

  Denote by ${\bold{a}}_{\Sigma}$ the two-sided ideal
of $\Sigma$ equal to the sum of nonprojective submodules   of
$\Soc \Sigma$, and let $\pi_{\Sigma}
:\Sigma \to \Sigma/{\bold{a}}_{\Sigma}$ be the natural projection.
By Proposition \ref{proj3} (a),
$\Sigma/{\bold{a}}_{\Sigma}$ is hereditary. 

\begin{prop}\label{quiver3} Let $G=(v(G),a(G))$ be the quiver 
given by $v(G) =v(H)$ and $ a(G) = a(H) - \{ \alpha_{d}\, |\, d \in
\mathcal{D} \}$, where $\alpha_{d}$'s are described in Lemma \ref{quiver2}(a). Let $\psi : kG \to \Sigma /{\bold{a}}_{\Sigma}$
be a unique morphism of $k$-algebras determined by $\psi( e_x) = \pi_{\Sigma}(\hat{e}_x)$ for $ x\in v(G)$, and
$\psi (\gamma) = \pi_{\Sigma}(\hat{\gamma})$ for $\gamma \in
a(G)$. Then $\psi$ is an isomorphism.
\end{prop}
\begin{pf} In view of Lemma \ref{quiver2},  $\psi$ is onto. To
show that $\psi$ is an isomorphism, it suffices to show that
$k[G]$ and $ \Sigma /{\bold{a}}_{\Sigma}$ have the same dimension
as $k$-spaces. The latter will be established if we show that $G$
is the quiver associated to the (hereditary) algebra $ \Sigma /{\bold{a}}_{\Sigma}$ (see
~\cite[Proposition III.1.13]{ars}).

Since $ \radic (\Sigma /{\bold{a}}_{\Sigma} ) =
{\bold{s}}/{\bold{a}}_{\Sigma}$, then $(\radic (\Sigma/
{\bold{a}}_{\Sigma}))^2 = ({\bold {s}} ^2 + {\bold{a}}_{\Sigma} )/
{\bold{a}}_{\Sigma}$ and we have 
\begin{equation}\label{quiver3.1}\Sigma /{\bold{a}}_{\Sigma} \biggr /({
\bold {s}}^2 + {\bold{a}}_{\Sigma} )/ {\bold{a}}_{\Sigma} \simeq
\Sigma /({\bold {s}}^2+ {\bold{a}}_{\Sigma} ) \simeq \Sigma /
{\bold {s}}^2 \biggr/({\bold {s}}^2+ {\bold{a}}_{\Sigma} )/{\bold
{s}}^2.
\end{equation} 
From Proposition \ref{proj3}(b) and Lemma \ref{quiver2}
it is clear that $\Sigma/{\bold {s}}^2 \biggr /({\bold {s}}^2+
{\bold{a}}_{\Sigma})/ {\bold {s}}^2 = \oplus k \tilde{\beta}$,
where $\beta \in a(H) - \{ \alpha_{d}\, |\, d \in \mathcal{D} \}$, and
$\tilde{\beta}$ is the image of $\hat{\beta}$ under the canonical
epimorphism 
$$ \varSigma /{\bold{a}}_{\varSigma} \to \dfrac{\varSigma/{\bold
{s}}^2 }{({\bold {s}}^2+ {\bold{a}}_{\varSigma})/ {\bold {s}}^2}
$$ 
(see Eq. (\ref{quiver3.1})). Therefore $\Sigma/{\bold {s}}^2 \biggr /({\bold {s}}^2+
{\bold{a}}_{\Sigma})/ {\bold {s}}^2$ has $G$ as its
quiver. It remains to use the fact that if $\Gamma$ is an
elementary algebra, then $\Gamma$ and $\Gamma / {\radic \Gamma}^2$
have the same quiver.
\end{pf}

\begin{lem}\label{quiver3a} In the setting of
Notation \ref{not1}, the triple $(G,{\mathcal D},\rho)$ satisfies the conditions
of Definition \ref{def1}.
\end{lem}
\begin{pf}  Viewing the isomorphism $\psi$ of Proposition \ref{quiver3} as identification, we may assume that
$\Sigma/{\bold{a}}_{\Sigma} = kG$. Since $\Sigma/{\bold{a}}_{\Sigma}$ is
finite-dimensional, $G$ has no oriented cycles. By Lemma \ref{quiver2} and
Proposition \ref{quiver3},  each $d \in \mathcal{D}$ is a sink in
$G$. If $y \in \Img {\rho}$ is a sink in $G$, then $S_y$ is a
projective $\Sigma/a_{\Sigma}$-module and nonprojective
$\Sigma$-module. In view of Proposition \ref{proj3}(c) and  equation (\ref{quiver1.1}), we get $y \in \mathcal{D}$.
\end{pf}

\begin{thm}\label{quiver4} In the  setting of Notation \ref{not1},
  $\Sigma = k[H, \sigma(H)] \simeq k(G,{\mathcal D},\rho)=k[H, \rho(H)]= \La$.
\end{thm}
\begin{pf} Recall that by Lemma \ref{quiver3a}
and Corollary \ref{quiver1}, $\La$ satisfies (F). Since $\langle\rho(H)\rangle \subset a_{\Sigma}$, Lemma
\ref{quiver2}(b) implies that there exists an epimorphism  $ \phi : \La \to \Sigma$ of $k$-algebras
satisfying $\phi\pi_{\rho} = \pi_{\sigma}$. To show $\phi$ is
an isomorphism, it suffices to show that so is the restriction $
\phi| \La {\bar{e}}_i= \phi_i : \La {\bar{e}}_i \to \Sigma
{\hat{e}}_i$, for each $i \in v(H)$. We proceed by induction on
$\dim_k\La \bar{e}_i$. If $\dim_k\La \bar{e}_i =1$, then
$\dim_k\Sigma \hat{e}_i =1$, and $\phi_i$ is an isomorphism. Suppose
$\dim_k\La \bar{e}_i= n >1$ and $\phi_j:\La \bar{e}_j \to
\Sigma \hat{e}_j$ is an isomorphism whenever $\dim_k\La \bar{e}_j
< n$. Since $\La \bar{e}_i$ is not simple, at least one arrow of
$H$ starts at $i$. Hence $\Sigma \, \hat{e}_i$ is not simple and ${\bold{s}}  \, \hat{e}_i \neq 0$.

Consider the exact commutative diagram of $\La$-modules

\begin{equation}\label{cd1}
\begin{CD}
0@>>> {\rad  \bar{e}_i}@>>>\La \bar{e}_i@>>>S_i@>>>0\\
@.@V{\phi_i|\rad  \bar{e}_i}VV@V{\phi_i}VV@|\\ 0@>>>{\bold {s}} \,
\hat{e}_i@>>>\Sigma \, \hat{e}_i@>>>S_i@>>>0
\end{CD}\ .
\end{equation}

\noindent Since $\Sigma$ satisfies (F), ${\bold {s}} \, \hat{e}_i$
is either projective or simple.

If ${\bold {s}}\,  \hat{e}_i$ is simple nonprojective,
$\ell( \Sigma \, \hat{e}_i)=2 $ and $i \in \mathcal{D}$.
By Proposition \ref{noloop}(c), $\dim_k\La \bar{e}_i=2
$, and $\phi_i$ is an isomorphism.

If ${\bold {s}} \,  {\hat{e}}_i$ is projective, $i
\notin \mathcal{D}$ and, by Proposition \ref{noloop}(d), $\rad
{\bar{e}}_i $ is projective. Let $\gamma_1 , \dots,\gamma_m $
be the arrows in $a(H)$ starting at $ i$ (remember, ${\bold{s}}  \, \hat{e}_i \neq 0$).  The  maps
$\overset{m}{\underset{u=1} {\oplus}} \La \bar{e}_{t(\gamma_u)}
\to \rad \bar{e}_i = \sum \limits_{u=1}^{m} \La \bar{\gamma}_u$
and $\overset{m}{\underset{u=1} {\oplus}} \Sigma \,
\hat{e}_{t(\gamma_u)} \to {\bold {s}} \,
  \hat{e}_i = \sum \limits_{u=1}^{m} \Sigma \,
\hat{\gamma}_u$ induced by the right multiplication by the
$\gamma_u$'s are projective covers, hence, isomorphisms, because
$\rad \bar{e}_i$ and  ${\bold {s}} \, \hat{e}_i$ are projective.
We obtain a commutative diagram
\begin{equation*}
\begin{CD}
{\overset{m}{\underset{u=1} {\oplus}} \La
\bar{e}_{t(\gamma_u)}}@>>>{\rad \bar{e}_i}\\
  @V{\overset{m}{\underset{u=1} {\oplus}}
\phi_{t(\gamma_u)}}VV @VV{\phi_i|\rad \bar{e}_i}V
\\ {\overset{m}{\underset{u=1} {\oplus}}
\Sigma \, \hat{e}_{t(\gamma_u)}}@>>>{\bold{s}}\,  \hat{e}_i
\end{CD} \ \ \ \ \ \ \ .
\end{equation*}
\noindent By the induction hypothesis, each $\phi_{t(\gamma_u)}$
is an isomorphism, hence so are $\overset{m}{\underset{u=1}
{\oplus}} \phi_{t(\gamma_u)}$ and $\phi_i | \rad  \bar{e}_i$.
Using commutative diagram (\ref{cd1}), we
see that $\phi_i$ is an isomorphism. \end{pf}

\begin{thm}\label{quiver5}
\begin{itemize}
\item [(a)] An elementary finite-dimensional $k$-algebra satisfies
(F) if and only if it is isomorphic to the algebra $k(G,{\mathcal D}, \rho)$
for some triple
  $(G,{\mathcal D}, \rho)$.
\item [(b)] For a given algebra satisfying (F), the triple $(G,{\mathcal D},
\rho)$ is unique up to isomorphism.
\end{itemize}
\end{thm}
\begin{pf} (a) This follows from
Corollary \ref{quiver1} and Theorem \ref{quiver4}.

(b) Let $\La = k[H, \rho(H)]$ and $\La' = k[H', \rho'(H')]$ be the
algebras associated with  triples $(G,{\mathcal D}, \rho)$ and $(G',{\mathcal D'}, \rho')$,
respectively. We show that if there exists an isomorphism of $k$-algebras $\psi:
\La \to \La'$, then there exists an isomorphism $\phi : G \to G'$
such that $\phi (\mathcal{D})= \mathcal{D'}$ and $\rho'\phi|{\mathcal D} =
\phi\rho$.

Since $\{ {\bar{e}}_x | x \in v(H) \}$ is a complete set of
primitive orthogonal idempotents in the  elementary algebra $\La$,
so is $ \{ \psi ({\bar{e}}_x) | x \in v(H) \}$ in $\La'$. It
follows that $\psi ({\bar{e}}_x) = {\hat{e}}_{x'} + r'$, for a
unique $x' \in v(H')$ and some $r' \in \radic \La'$, whence the map $\phi
: v(H) \to v(H')$ given by $\phi(x) =x'$ is a bijection. In view
of ~\cite[Proposition 1.14]{ars}, which is also valid for an
elementary $k$-algebra, $\phi$ can be extended to
an isomorphism $\phi: H \to H'$. Since $v(G) = v(H)$ and $a(H) =
a(G) \cup \{ \alpha_{d} \,|\, d \in \mathcal{D} \}$, where
$\alpha_{d}$ is the unique arrow in $H$ starting at $d$, and
since $v(G'), v(H'), a(G')$, and $a(H')$ are similarly related, it
suffices to show that $\phi|\mathcal{D}: \mathcal{D}\to
\mathcal{D'}$ is a bijection. Indeed, assuming the map is a bijection, note that since $\phi$ is an
isomorphism of quivers, $\phi(\alpha_{d})$ is the unique arrow in
$a(H')$ starting at $\phi(d)$, so we must have $\phi(\alpha_{d}) =
\alpha_{\phi(d)}$. Hence $\phi| a(G): a(G) \to a(G')$ is a
bijection,  and $\rho'\phi|{\mathcal D} =
\phi\rho$, whence $\phi: (G, {\mathcal D},\rho) \to (G', {\mathcal D'},\rho')$ is an isomorphism.

To show that $\phi|\mathcal{D}: \mathcal{D}\to
\mathcal{D'}$ is a bijection, we prove that $\phi(\mathcal{D}) \subset \mathcal{D'}$; the
proof of $\phi^{-1}(\mathcal{D'}) \subset \mathcal{D}$ is similar.
If $ d \in \mathcal{D}$, Proposition \ref{noloop} (c) implies that
$\La {\bar{e}}_{d}$ is an indecomposable projective module of
length 2 with a nonprojective simple socle, and so is $ \psi(\La
{\bar{e}}_{d})= \La' \psi({\bar{e}}_d) = \La'({\hat{e}}_{\phi(d)}
+ r')=\La' {\hat{e}}_{\phi(d)} $, where $r' \in \radic \La'$.
Using parts (c) and (d) of  Proposition \ref{noloop}, we get
$\phi(d) \in \mathcal{D'}$.
\end{pf}

Having shown that a finite-dimensional elementary $k$-algebra
satisfies (F) if and only if it is isomorphic to a quadratic
monomial algebra $k(G,{\mathcal D},\rho)$, we explain how to verify
whether a given quadratic monomial algebra satisfies (F). Our
description is similar to that in ~\cite{br} for the algebras stably
equivalent to hereditary.

\begin{nt}\label{not2} Let $F$ be a
finite quiver  with a set of quadratic monomial relations
$\sigma(F)$, and denote by $\hat{p}$ the image of a path $p$ in $F$
under the natural projection $kF \to k[F, \sigma(F)]$. We set 
\newline\centerline{$\mathcal{X}$ =$ \{ x \in v(F)\, |\, \beta \alpha \in \sigma(F)$ for
some $\alpha, \beta \in a(F)$ with $ s(\alpha) =x \}$.}
\end{nt}

\begin{prop}\label{quiver6} In the setting of Notation \ref{not2},
$\Sigma = k[F, \sigma(F)] \simeq k(G,{\mathcal D},
 \rho)$ for some triple $(G,{\mathcal D},
\rho)$  if and only if the following conditions hold.
\begin{itemize}
   \item [(a)] For all $ x \in \mathcal{X}$, there is precisely
   one $\alpha \in a(F)$ satisfying $s(\alpha) =x$.
   \item [(b)] If $\alpha \in a(F)$ satisfies $s(\alpha) \in
   \mathcal{X}$, then  for all $\beta \in a(F)$ satisfying $s(\beta) =
   t(\alpha)$, we have $\beta\alpha \in \sigma(F)$.
   \item [(c)] Each oriented cycle passes through a vertex in
   $\mathcal{X}$.
\end{itemize}
\end{prop}
\begin{pf} Suppose $\Sigma= k[F, \sigma(F)] \simeq
k(G, {\mathcal D},
\rho)$.  By Theorem
\ref{quiver5}(a), $\Sigma$ satisfies (F). Since $\Sigma$ is finite-dimensional, $\mathcal{X}$ must satisfy
(c). If $x \in \mathcal{X}$, there exist $\alpha,
\beta \in a(F)$ with $s(\alpha) =x$  and $\beta\alpha \in
\sigma(F)$. Since $\Sigma$ is a monomial algebra, $\Sigma \,
\hat{\alpha}$ is a direct summand of ${\bold{s}} \,
{\hat{e}}_{x}$ where ${\bold{s}}=\radic\Sigma$, and $\Sigma \, \hat{\alpha}$ is not projective
because $\hat{\beta} \hat{\alpha} =0$. Since $\Sigma$ satisfies
(F), $\Sigma \, \hat{\alpha} \simeq {\bold{s}} \, {\hat{e}}_{x}$ is
simple, so (a) and (b) follow.

Suppose $F$ and $ \sigma(F)$ satisfy (a), (b), and
(c). Denote by $\alpha_{x}$ the unique arrow with
$s(\alpha_{x})=x \in \mathcal{X}$. Let $G$ be the quiver
defined by $v(G) =v(F)$ and $ a(G) = a(F) - \{ \alpha_x \,|\, x \in
\mathcal{X} \}$. By (c), $G$ has no
oriented cycles. Define $\rho: \mathcal{X}$ $\to v(G)$ by $\rho(x)
= t(\alpha_x)$. By (a), $\mathcal{X}$
consists of sinks in $G$. Suppose $y \in \Img \rho$ is a sink in
$G$. By (b), there is $\beta \in a(F)$ satisfying $s(\beta)
=y$ and, by assumption, $\beta \notin a(G)$. Therefore $y \in
\mathcal{X}$.  We have just proved that $(G, {\mathcal D},\rho)$ satisfies the
conditions of Definition \ref{def1}. It remains to use
(b) to conclude that $k(G,{\mathcal D},\rho) = k[F, \sigma(F)]$.
\end{pf}

We now determine when both $\La$ and $\La^{op}$ satisfy (F). In the sequel, $H^{op}$ stands for the opposite quiver of $H$, and if $\gamma:x\to y$ is an arrow of $H$, we  denote by $\gamma^{op}:y\to x$ the corresponding arrow of  $H^{op}$.

\begin{prop}\label{quiver7} An elementary finite-dimensional
$k$-algebra $\La$ and its opposite satisfy (F) if and only if $
\La \simeq k(G,{\mathcal D}, \rho)$ and the following
conditions hold.
\begin{itemize}
   \item [(a)] Each $ y \in \Img {\rho}$ is a source in $G$.
   \item [(b)] If $\beta \in a(G)$ with  $s(\beta)
\in \Img {\rho}$ then $\beta$ is the only arrow of $G$
    that ends at $t(\beta)$.
    \item[(c)] If $\rho(d_1)=\rho(d_2)$  and $d_1\ne d_2$, then $d_i\not\in  \Img {\rho},\ i=1,2$.
\end{itemize}
\end{prop}
\begin{pf} By Theorem \ref{quiver5}(a), it suffices to
show that if $\La \simeq k(G,{\mathcal D}, \rho)=k[H,\rho(H)]$, 
then $\La^{op}$ satisfies (F) if and only if (a), (b) and (c) hold. In view of
Definition \ref{def2}, $\La^{op} \simeq k[F, \sigma(F)]$ where $F=H^{op}$ and $\sigma(F) = \{ \alpha^{op}\beta^{op} \,|\, \beta
\alpha \in \rho(H) \}$. Using Notation \ref{not2}, we
have $\mathcal{X}$ = $\{ t(\beta)\,|\, \beta\alpha
\in \rho(H)\}$.

Suppose $\La^{op}$ satisfies (F). If $\beta, \beta' \in a(G)$ with
$t(\beta) = t(\beta')$ and $s(\beta) \in \Img {\rho}$, then $
\beta\alpha \in \rho(H)$ for some $\alpha \in a(H)$ so that
$s(\beta^{op}) = s({\beta'}^{op}) \in \mathcal{X}$. Applying Proposition
\ref{quiver6}(a) to $\La^{op}$, we get $\beta^{op} =
{\beta'}^{op}$ whence $\beta = \beta'$. Thus (b) holds. Let
$\gamma \in a(G)$ satisfy $t(\gamma) \in \Img {\rho} $. Then
$t(\gamma) = t(\alpha)$ where $\alpha \in a(H)$ so $\beta\alpha
\in \rho(H)$ for some $\beta \in a(H)$. Hence $s(\beta^{op})
=t(\beta) \in \mathcal{X}$ and $s(\gamma^{op}) =s(\alpha^{op}) =
t(\beta^{op})$. Applying Proposition \ref{quiver6}(b) to
$\La^{op}$, we get $\gamma^{op}\beta^{op} \in \sigma(F)$ whence $\beta\gamma \in \rho(H)$. The latter contradicts
Theorem \ref{quiver5}(b) because $\gamma \in a(G)$, so (a) holds.
Let $d_1, d_2\in \mathcal{D}$ be distinct vertices satisfying $\rho(d_1)=\rho(d_2)=x$.  If, say,  $d_1=\rho(d)$ for some $d\in \mathcal{D}$, then by Definition \ref{def2} there exist arrows $\alpha_{d_1} : d_1\to x,\ \alpha_{d_2}: d_2 \to x$, and $\alpha_d: d \to d_1$  in $a(H)$, and we have $\alpha_{d_1}\alpha_d \in
\rho(H)$ whence $\alpha_d^{op}\alpha_{d_1}^{op} \in
\sigma(F)$ and $x\in\mathcal X$.   This contradicts Proposition \ref{quiver6}(a): $\alpha_{d_1}^{op}, \alpha_{d_2}^{op}\in a(F) $ are distinct and $s(\alpha_{d_1}^{op})=s(\alpha_{d_2}^{op})=x$. We conclude that (c)
holds.

Suppose (a), (b) and (c) hold. To show
$\La^{op}$ satisfies (F), we verify the conditions of Proposition
\ref{quiver6}. Since $\La$ satisfies (F), Proposition \ref{noloop}(a) implies that each oriented cycle in $H$ passes through a vertex
in $\mathcal{D}$. By Lemma \ref{quiver2}(a), for each $d \in
\mathcal{D}$, there is a unique arrow $\alpha_{d} \in a(H)$ with $s(\alpha_{d}) =d$ and $t(\alpha_{d}) \in \Img {\rho}$. Hence
each oriented cycle  contains an arrow $\beta$ with $s(\beta) =
t(\alpha_d) \in \Img {\rho}$. By Definition \ref{def2}, $\beta
\alpha_d \in \rho(H)$ whence $t(\beta) \in \mathcal{X}$. Thus
each oriented cycle in $H$ passes through a vertex in
$\mathcal{X}$, and so does each cycle in $H^{op}$. Hence
condition (c) of Proposition \ref{quiver6} holds.

Let $x\in\mathcal X$, then $x=t(\beta)$ where $\beta\alpha\in\rho(H)$ for some $\alpha,\beta\in a(H)$ and $s(\beta)=t(\alpha)\in\Img\rho$.  Let $\gamma^{op}\in a(H^{op})$ with $s(\gamma^{op})=t(\gamma)=x$.  To show $\gamma=\beta$ and, hence, $\gamma^{op}=\beta^{op}$, assume, to the contrary, that $\gamma\not=\beta$.  If $\beta\in a(G)$ then, by (b), $\gamma\not\in a(G)$ whence $x=\rho(s(\gamma))$.  By (a), $x$ is a source in $G$, a contradiction.  If $\beta\not\in a(G)$ then
$x=\rho(s(\beta))$ is a source in $G$, so $\gamma\not\in a(G)$.  By Definition \ref{def2},  $s(\beta)\ne s(\gamma)$ which, together with $x=\rho(s(\beta))=\rho(s(\gamma))$ and $s(\beta)\in\Img\rho$, contradicts (c).  Thus we must have $\gamma=\beta$, so condition (a) of Proposition \ref{quiver6} holds.

Suppose $\beta^{op}, \gamma^{op} \in a(H^{op})$ satisfy
$s(\beta^{op}) \in \mathcal{X}$ and $s(\gamma^{op}) =
t(\beta^{op})$. Then $t(\beta) \in \mathcal{X}$ and there exists
$\alpha \in a(H)$ satisfying $\beta \alpha \in \rho(H)$. It
follows that $t(\gamma) = s(\beta) \in \Img {\rho}$ so that, by
(a), $\gamma \notin a(G)$. By Definition \ref{def2}, $\beta \gamma
\in \rho(H)$ whence $\gamma^{op} \beta^{op} \in \sigma (F)$.  Thus condition (b) of Proposition 3.8 holds.
\end{pf}

We finish this section with examples of elementary $k$-algebras satisfying (F).  In view of Theorem \ref{quiver5}, it suffices to indicate the triple $(G,{\mathcal D}, \rho)$ that determines such an algebra.

\begin{ex}\label{quiver8}  (a) Let $(G,{\mathcal D}, \rho)$ be the triple given by $v(G)=\{1\}$, $a(G)=\emptyset$, ${\mathcal D}=v(G)$, and $\rho(1)=1$.  Then $\La_1=k(G,{\mathcal D}, \rho)=k[H, \rho(H)]$ is the path algebra of the quiver $H$ with $v(H)=v(G)$ and $a(H)=\{\alpha_1:1\to 1\}$ modulo the ideal generated by $ \rho(H)=\{\alpha_1^2\}$.  Since $\La_1$ is commutative, its opposite algebra satisfies (F).\vskip.05in

(b) Let  $v(G)=\{1,2\}$, $a(G)=\emptyset$, ${\mathcal D}=v(G)$, and $\rho(1)=\rho(2)=2$.  Then $v(H)=v(G)$, $a(H)=\{\alpha_1:1\to 2,\ \alpha_2:2\to 2\}$, and $ \rho(H)=\{\alpha_2\alpha_1,\alpha_2^2\}$.  Since the triple $(G,{\mathcal D}, \rho)$ does not satisfy condition (c) of Proposition \ref{quiver7}, the opposite algebra of $\La_2=k(G,{\mathcal D}, \rho)$ does not satisfy (F).\vskip.05in

(c) Let $v(G)=\{1,2\}$, $a(G)=\{\beta:1\to 2\}$, ${\mathcal D}=\{2\}$, and $\rho(2)=2$.  Then $v(H)=v(G)$, $a(H)=\{\beta:1\to 2,\ \alpha_2:2\to 2\}$, and $ \rho(H)=\{\alpha_2^2\}$.  Since the triple $(G,{\mathcal D}, \rho)$ does not satisfy condition (a) of Proposition \ref{quiver7}, the opposite algebra of $\La_3=k(G,{\mathcal D}, \rho)$ does not satisfy (F).\vskip.05in

(d) Suppose that $G$ and  ${\mathcal D}$ are the same as in (c) but $\rho(2)=1$.  Then $v(H)=v(G)$, $a(H)=\{\beta:1\to 2,\ \alpha_2:2\to 1\}$, and $ \rho(H)=\{\beta\alpha_2\}$.  By Proposition \ref{quiver7}, the opposite algebra of $\La_4=k(G,{\mathcal D}, \rho)$ satisfies (F).\vskip.05in

Since $\gd\La_4<\infty$, the support variety of each module in ${\mathrm{mod}}\La_4$ is trivial.  The algebras $\La_1,\ \La_2$, and $\La_3$ have a lot in common: all three are of infinite global dimension, and each module in ${\mathrm{mod}}\La_i,\ i=1,2,3$, is eventually periodic by 
Corollary \ref{proj1.5}.  All three are quadratic monomial, hence, Koszul algebras so the cohomology ring ${\mathrm HH}^*(\La_i)$ maps onto the graded center of the Ext-algebra $E(\La_i)$.  Since the graded center of $E(\La_2)$ is the ground field $k$, the support variety of each module in  ${\mathrm{mod}}\La_2$ is trivial.  However, when $i$ equals  $1$ or $3$, the graded center is the polynomial algebra in one variable with a suitable grading, and  $E(\La_i)$ is a finitely generated module over the graded center.  In this setting the methods of  \cite{snsol04} apply, and using these methods, it is possible to show that the unique simple $\La_1$-module and the simple $\La_3$-module associated to the vertex $2$ have nontrivial support varieties  ~\cite{s}.
\end{ex} 

\section{The Ext-algebra of $\La$ satisfying (F)}
In this section, $\La$ is an elementary finite-dimensional $k$-algebra satisfying (F). We study the Ext-algebra $E(\La)={\overset{\infty}{\underset{i=0}{\oplus}}}\Ext^{i}_{\La}(\La/\rad,\La/\rad)$.

\begin{thm}\label{ext1}
Let $\La \simeq k(G,{\mathcal D}, \rho) =k[H, \rho(H)]$ as
described in Definitions \ref{def1} and \ref{def2}.
\begin{itemize}
\item[(a)] $ E(\La) \simeq k[H, \eta(H)]$ where $\eta (H) = \{
\beta \gamma \, |\, \beta  \in a(H), \gamma \in a(G) \}$. 
\item[(b)] $E(\La)$ is left
noetherian.
\end{itemize}
\end{thm}
\begin{pf} (a) The underlying quiver and relations of the Ext-algebra of a quadratic monomial algebra are described in ~\cite[Proposition 2.2]{gz}.  Applying this description to $\La$, we obtain the statement as an immediate consequence of Theorem \ref{quiver5}.

(b) Let $Q$ be
the quiver with $v(Q) = v(H)$ and $a(Q) = \{
\alpha_{d}\,| \,d \in \mathcal{D} \}$.  Let $\psi :kH \to kQ$ be a  unique epimorphism of $k$-algebras determined by $\psi(e_i) =e_i$ for all $i \in v(H)$;  $\psi(\beta) =
0$ for all $\beta \in a(G)$; and $\psi(\alpha_{d}) =\alpha_{d}$ for all
$d \in \mathcal{D}$. Since $\langle\eta(H)\rangle
\subset \Ker \psi$, there exists a unique epimorphism of algebras
$\phi:k[H, \eta(H)] \to kQ$ satisfying $\psi=\phi\pi$ where $\pi:kH\to k[H, \eta(H)]$ is the natural projection. To prove that
$E(\La)$ is left noetherian, we show that $kQ$ is left noetherian, and $\Ker \phi$ is
finite-dimensional over $k$.

For $kQ$, we use ~\cite[Theorem
8]{ls}, which says that the path algebra $k\Gamma$ of a finite
quiver $\Gamma$ is left noetherian if and only if $\Gamma$
satisfies the following condition. If there is an oriented cycle
passing through a vertex $i$ of $\Gamma$, then only one arrow
starts at $i$.  By Definition \ref{def2}, $Q$
satisfies this condition.

It is easy to see that, first, $\Ker \phi$ is
spanned over $k$ by the elements $\pi(p)$ where  $p$ is a path in $\Ker
\psi$, and, second, $p\in\Ker\psi$  if and only if $p$ contains an arrow of $G$. Suppose $p \in \Ker \psi$ is an oriented cycle in $H$. Since $G$ has
no oriented cycles, $p$ must contain an arrow $\alpha_{d},d\in \mathcal{D}$, and, hence,  a subpath of the form
$\alpha_{d} \gamma$ with $\gamma\in a(G)$, which implies $\pi(p) =0$. Therefore $\Ker
\phi$ is spanned by the elements $\pi(p)$ where $p\in \Ker \psi$ and no subpath of $p$ is an oriented cycle.  Since there are only finitely many such paths $p$, $\Ker \phi$ is finite-dimensional.
\end{pf}

We note that the Ext-algebra $E(\La)$ need not be right noetherian.  For example, for the algebra $\La_2$ of Example \ref{quiver8}(b), $E(\La_2)$ is the path algebra of the quiver with the set of vertices $\{1,2\}$ and  the set of arrows $\{\alpha_1:1\to 2,\ \alpha_2:2\to 2\}$.  It is not right noetherian because  the loop $\alpha_2$ passes through the vertex $2$, and more than one arrow ends at $2$.


\begin{thebibliography}{F7}

\bibitem[1]{ar}
M. Auslander and I. Reiten, Stable equivalence of artin algebras,
{\it in} \lq\lq Proceedings of the Conference on Orders, Group
Rings and Related Topics," Ohio State University, 1972, Lecture
Notes in Mathematics, Vol. 353, Springer-Verlag,
Berlin-Heidelberg-New York, 1973, pp. 8-64.

\bibitem[2]{ars}
M. Auslander, I. Reiten and S. O. Smal\o, \lq\lq Representation
Theory of Artin Algebras," Cambridge Studies in Advanced
Mathematics, Vol. 36, Cambridge University Press, New York, 1994.

\bibitem[3]{av89a} L. L. Avramov, Homological
asymptotics of modules over local rings, {\it in} "Commutative
algebra (Berkeley, CA, 1987)," Math. Sci. Res. Inst. Publ., {\bf
15}, Springer, New York, 1989, pp. 33-62.

\bibitem[4]{be91}
D. J. Benson, \lq\lq Representations and Cohomology II," Cambridge
Studies in Advanced Mathematics, Vol. 31, Cambridge University
Press, Cambridge, 1991.

\bibitem[5]{br} K. Bongartz and C. Riedtmann, Alg\`ebres
stablement h\'er\'editaires. (French) {\it C.R. Acad. Sci. Paris S\'er. A-B} {\bf 288} (1979), no. 15, 703--706.

\bibitem[6]{ehsst} K. Erdmann, M. Holloway,
N. Snashall, \O. Solberg, R. Taillefer, Support varieties for
selfinjective algebras,  {\it $K$-Theory} {\bf 33} (2004), 67-87.


\bibitem[7]{gz} E.L. Green and D. Zacharia, The Cohomology Ring of a
Monomial Algebra, {\it manuscripta math.} 85, 11-23 (1994).

\bibitem[8]{jk1}
S. Jagadeeshan and M. Kleiner, Stable artin algebras: the
transpose and the global dimension, {\it in} \lq\lq Seventh
Conference on Representations of Algebras (ICRA VII)," Canadian
Mathematical Society Conference Proceedings, Vol. {\bf 18}, 1996,
pp. 343-351.


\bibitem[9]{jk2}
S. Jagadeeshan and M. Kleiner, Structure of projectively stable
artinian rings, {\it J. Algebra} {\bf 197} (1997), 92-126.

\bibitem[10]{jk3}
S. Jagadeeshan and M. Kleiner, Projectively stable artin algebras:
representation theory and quivers with relations, {\it Comm.
Algebra} {\bf 27(5)} (1999), 2277-2316.

\bibitem[11]{ls}
S. Liu, Isomorphism problem for tensor algebras over valued graphs,
{\it Sci.China Ser.A}, {\bf 34}(1991), no. 3, pp. 267-272.

\bibitem[12]{mz}
R. Mart\'{i}nez-Villa and D. Zacharia, Approximations with modules having linear resolutions.  {\it J. Algebra} {\bf 266} (2003), no. 2, 671--697.

\bibitem[13]{pl}
M. I. Platzeck, Representation theory of algebras stably
equivalent to an hereditary artin algebra, {\it Trans. Amer. Math.
Soc} {\bf 238} (1978), 89-128.

\bibitem[14]{snsol04} N. Snashall and \O. Solberg,
Support varieties and Hochschild cohomology rings,  {\it Proc.
London Math. Soc. (3)} {\bf 88} (2004), no. 3, 705--732.

\bibitem[15]{s} \O. Solberg,
private communication.


\end{thebibliography}
\end{document}